\documentclass[10.5pt]{amsart}
\usepackage[a4paper,margin=1.00in]{geometry}
\usepackage{graphicx}
\usepackage[colorlinks, citecolor=blue, linkcolor=black]{hyperref}

\title[A one-sentence proof of the Bernstein type theorem for special Lagrangian equation in two dimensions]{A one-sentence proof of the Bernstein type theorem \\ for special Lagrangian equation in two dimensions}
\author{Hojoo Lee}
\email{momentmaplee@gamil.com}

\newtheorem{thm}{Theorem}

\theoremstyle{definition}

\numberwithin{equation}{section}

\begin{document}

\maketitle

The aim of this article is to reveal the equivalence of J\"{o}rgens' Theorem  \cite{Jorgens 1954} and Fu's Theorem \cite{Fu 1998}.
  
 \begin{thm} [\textbf{J\"{o}rgens \cite{Jorgens 1954}}]\label{jorgens} The only entire ${\mathcal{C}}^{2}$ solutions of the Monge-Amp\`{e}re equation
\[
  h_{xx} h_{yy} -{h_{xy}}^{2}=1, \quad (x, y) \in {\mathbb{R}}^{2}
\]
are quadratic polynomial functions. 
 \end{thm}
 
In  1957, J. C. C. Nitsche \cite{Nitsche 1957} combined the Heinz observation (\cite[p. 133]{Jorgens 1954},  \cite[Lemma 4.4]{Osserman 1986}) and J\"{o}rgens' Theorem to present a single page proof of Bernstein's beautiful theorem that the only entire minimal graphs in  ${\mathbb{R}}^{3}$ are planes.  L. Fu \cite{Fu 1998} exploited Osserman's powerful Bernstein type result (\cite{Osserman 1966}, \cite[Theorem 5.1]{Osserman 1986}) 
for two dimensional entire minimal graphs in  ${\mathbb{R}}^{N \geq 3}$ to characterize entire solutions of the special Lagrangian equation 
 \cite{Harvey Lawson 1982, Warren Yuan 2009}, which combines Laplace equation and unimodular Hessian equation.
  
 \begin{thm}[\textbf{Fu \cite{Fu 1998}}] \label{fu} Let $\theta \in \mathbb{R}$ be a constant. If an entire function ${\mathcal{F}}$ of 
 class ${\mathcal{C}}^{2}$ satisfies the SLE
\[
  \cos \theta \left[ \, {\mathcal{F}}_{xx} + {\mathcal{F}}_{yy}  \, \right] +  \sin \theta  \left[ \,  {\mathcal{F}}_{xx} {\mathcal{F}}_{yy} -{{\mathcal{F}}_{xy}}^{2} - 1 \, \right] =0, \quad (x, y) \in {\mathbb{R}}^{2}, 
\]
then the function ${\mathcal{F}}$ is harmonic or quadratic.
 \end{thm}
 
The linearization in the proof \cite{Jorgens 1954} of J\"{o}rgens' Theorem  via the partial Legendre transformation can be viewed
as a geometrical fact that the gradient graph of any solution of the special Lagrangain equation becomes a Lagrangian minimal graph 
 in ${\mathbb{R}}^{4}$, which can be identified as a holomorphic curve in ${\mathbb{C}}^{2}$.  
 
 \begin{proof}[\textbf{Alternative Proof of Theorem \ref{fu}.}]
As in \cite[Lemma 18]{Lee 2013}, when $\sin \theta \neq 0$, observe that the function
 \[ 
 h(x,y) = \left( \cos \theta \right) \frac{x^{2}+y^{2}}{2} + \left( \sin \theta \right) {\mathcal{F}}(x,y) 
 \]
solves the equation
$h_{xx} h_{yy} -{h_{xy}}^{2} = 1 + \sin \theta  \left\{ \,  \cos \theta \left[   {\mathcal{F}}_{xx} + {\mathcal{F}}_{yy}   \right] +  \sin \theta  \left[    {\mathcal{F}}_{xx} {\mathcal{F}}_{yy} -{{\mathcal{F}}_{xy}}^{2} - 1   \right]  \, \right\} = 1$.
 \end{proof}
 
 The equivalence of Theorems \ref{jorgens} and \ref{fu} seems to have gone unnoticed in the literature and may have some potential 
 applications to more general (Hessian type) equations. For several modern Bernstein type results for higher dimensional special 
 Lagrangian graphs, we refer readers to \cite{Jost Xin 2002, Tsui Wang 2002, Yuan 2002}. For more recent applications of Osserman's 
 higher codimensional generalization (\cite{Osserman 1966}, \cite[Theorem 5.1]{Osserman 1986}) of Bernstein's Theorem, we also refer 
 interested readers to \cite[Theorem 1.1 and Corollary 1.4]{H SH V 2009} and \cite[Theorem 5.2]{Jost Xin Yang 2002}.

\bigskip

\end{document}